\documentclass[11pt]{article}
\usepackage{amsmath}
\usepackage{amsfonts}
\usepackage{amssymb}
\usepackage{enumerate}
\usepackage{amsthm}
\usepackage{authblk}
\usepackage{graphicx}
\usepackage{verbatim}
\newtheorem{myproposition}{Proposition}[section]
\newtheorem{mytheorem}[myproposition]{Theorem}
\newtheorem{mylemma}[myproposition]{Lemma}

\newtheorem{mycorollary}[myproposition]{Corollary}

\newtheorem{myobservation}[myproposition]{Observation}

\def\gr{\mathcal{G}}

\makeatletter
\def\imod#1{\allowbreak\mkern10mu({\operator@font mod}\,\,#1)}
\makeatother
\begin {document}

\def\zet{\mathbb{Z}}

\title{Group Irregular Labelings of Disconnected Graphs}


\author[1,3]{Marcin Anholcer}
\author[2,3]{Sylwia Cichacz}
\affil[1]{\scriptsize{}Pozna\'n University of Economics, Faculty of Informatics and Electronic Economy}
\affil[ ]{Al. Niepodleg{\l}o\'sci 10, 61-875 Pozna\'n, Poland, \textit{m.anholcer@ue.poznan.pl}}
\affil[ ]{}
\affil[2]{AGH University of Science and Technology, Faculty of Applied Mathematics}
\affil[ ]{Al. Mickiewicza 30, 30-059 Krak\'ow, Poland, \textit{cichacz@agh.edu.pl}}
\affil[ ]{}
\affil[3]{University of Primorska, Faculty of Mathematics, Natural Sciences and Information Technologies}
\affil[ ]{Glagolja\v{s}ka 8, SI-6000 Koper, Slovenia, \{marcin.anholcer, sylwia.cichacz-przenioslo\}@famnit.upr.si}

\maketitle

\begin{abstract}
We investigate the \textit{group irregularity strength} ($s_g(G)$) of graphs, i.e. the smallest value of $s$ such that taking any Abelian group $\gr$ of order $s$, there exists a function $f:E(G)\rightarrow \gr$ such that the sums of edge labels at every vertex are distinct. We give the exact values and bounds on $s_g(G)$ for chosen families of disconnected graphs. In addition we present some results for the \textit{modular edge gracefulness} $k(G)$, i.e. the smallest value of $s$ such that there exists a function $f:E(G)\rightarrow \zet_s$ such that the sums of edge labels at every vertex are distinct.
\end{abstract}

\noindent\textbf{Keywords:} irregularity strength, graph weighting,
graph labeling, Abelian group\\
\noindent\textbf{MSC:} 05C15, 05C78.

\section{Introduction}

It is easy exercise to show that in any simple graph $G$ there are at least two vertices of the same degree. The situation changes if we consider an edge labeling $f:E(G)\rightarrow \{1,\dots,s\}$ and calculate weighted degree (or weight) of each vertex $x$ as the sum of labels of all the edges incident to $x$. The labeling $f$ is called \textit{irregular} if the weighted degrees of all the vertices are distinct. The smallest value of $s$ that allows some irregular labeling is called \textit{irregularity strength of $G$} and denoted by $s(G)$.

The problem of finding $s(G)$ was introduced by Chartrand et al. in \cite{ref_ChaJacLehOelRuiSab1} and investigated by numerous authors. Best published result due to Kalkowski et al. (see \cite{ref_KalKarPfe1}) is $s(G)\leq 6n/\delta$. There are some signals that it was recently improved by Przyby{\l}o (\cite{ref_Prz3}) for dense graphs of sufficiently big order ($s(G)\leq 4n/\delta$ in this case). Exact value of $s(T)$ for a tree $T$ was investigated e.g. by Aigner and Triesch (\cite{ref_AigTri}), Amar and Togni (\cite{ref_AmaTog}), Ferrara et al. (\cite{ref_FerGouKarPfe}) and Togni (\cite{ref_Tog1}).

Lo in \cite{ref_Lo} defined \textit{edge graceful labelings}. A graph $G$ of order $n$ and size $m$ is edge-graceful if there exists a bijective mapping  $f: E(G) \rightarrow\{1,2,\dots,m\}$ such that the induced vertex labeling
$f^\prime : V (G)\rightarrow \{0, 1,2,\ldots, n-1\}$ defined as 

$$f^\prime(u) =\sum_{v\in N(u)}f(uv)\mod n$$
is a bijection. In the above definition of an edge-graceful labeling of a connected graph $G$ of order $n>2$ and size $m$, the edge labeling $f$ is required to be a bijection. However, since the induced vertex labels $f'(u)$ are obtained by summation in $\zet_n$, the mapping $f$ is actually a mapping from $E(G)$ to $\zet_n$ and is in general not one-to-one. It is easy to see that the mapping  $f:E(G) \rightarrow \zet_n$ is one-to-one if and only if $m = n-1$ or $m = n$ (\cite{ref_Jon}). Jones combined the concepts of graceful labeling and modular edge coloring into labeling called  a \emph{modular edge-graceful labeling} (\cite{ref_Jon,ref_JonKolOkaZha2, ref_JonZha}). He defined the \textit{modular edge-gracefulness of graphs} as the smallest integer $k(G)=k \geq n$ for which there exists an edge labeling $f:E(G)\rightarrow \zet_k$ such that the induced vertex labeling $f^\prime : V (G)\rightarrow \zet_k\}$ defined by
$$
f^\prime(u) =\sum_{v\in N(u)}f(uv)\mod k
$$
is one-to-one.

Assume $\gr$ is an Abelian group  of order $m\geq n$ with the operation denoted by $+$ and identity element $0$. For convenience we will write $ka$ to denote $a+a+\dots+a$ (where element $a$ appears $k$ times), $-a$ to denote the inverse of $a$ and we will use $a-b$ instead of $a+(-b)$. Moreover, the notation $\sum_{a\in S}{a}$ will be used as a short form for $a_1+a_2+a_3+\dots$, where $a_1, a_2, a_3, \dots$ are all the elements of the set $S$.
The order of the element $a\neq 0$ is the smallest $r$ such that $ra=0$. It is well-known by Lagrange Theorem that $r$ divides $|\gr|$ \cite{ref_Gal}.

We define edge labeling $f:E(G)\rightarrow \gr$ leading us to the weighted degrees defined as the sums (in $\gr$):
$$
w(v)=\sum_{v\in N(u)}f(uv)
$$

The concept of $\gr$-irregular labeling is a generalization of a modular edge-graceful labeling. In both cases the labeling $f$ is called \textit{$\gr$-irregular} if all the weighted degrees are distinct. However, the \textit{group irregularity strength} of $G$, denoted $s_g(G)$, is the smallest integer $s$ such that for every Abelian group $\gr$ of order $s$ there exists $\gr$-irregular labeling $f$ of $G$. Thus the following observation is true.
\begin{myobservation}
For every graph $G$ with no component of order less than $3$, $k(G)\leq s_g(G)$.

\end{myobservation}
The following theorem, determining the value of $s_g(G)$ for every connected graph $G$ of order $n\geq 3$, was proved by Anholcer and Cichacz \cite{ref_AnhCic1}.

\begin{mytheorem}[\cite{ref_AnhCic1}]\label{AnhCic1}
Let $G$ be an arbitrary connected graph of order $n\geq 3$. Then
$$
s_g(G)=\begin{cases}
n+2,&\text{if   } G\cong K_{1,3^{2q+1}-2} \text{  for some integer   }q\geq 1\\
n+1,&\text{if   } n\equiv 2 \imod 4 \wedge G\not\cong K_{1,3^{2q+1}-2} \text{  for any integer   }q\geq 1\\
n,&\text{otherwise}
\end{cases}
$$
\end{mytheorem}

In \cite{ref_Jon} it was proved in turn that for every connected graph $G$ of order $n\geq 3$ 
$$
k(G)=\left\{\begin{array}{lll}
n,&\text{if}& n \not \equiv 2 \imod 4, \\
n+1,&\text{if}& n\equiv 2 \imod 4.
\end{array}
\right.
$$

In order to distinguish $n$ vertices in an arbitrary (not necessarily connected) graph we need at least $n$ distinct
elements of $\gr$. However, $n$ elements are not always enough, as shows the
following lemma.

\begin{mylemma}[\cite{ref_AnhCic1}]\label{lemma_below}
Let $G$ be of order $n$, if $n \equiv 2 \imod 4$, then $s_g(G) \geq n+1$.
\end{mylemma}

In this paper we generalize the above results for some families of disconnected graphs. Our main results are three following theorems.

\begin{mytheorem}\label{main}
Let $G$ be a graph of order $n$ with no component of order less than $3$ and with all the bipartite components having both color classes of even order. Then:
$$
\begin{array}{lll}
s_g(G)=n,& \text{if} &n \equiv 1\imod 2,\\
s_g(G)=n+1,& \text{if} &n\equiv 2 \imod 4,\\
s_g(G)\leq n+1, & \text{if} &n\equiv 0 \imod 4.
\end{array}
$$
\end{mytheorem}

\begin{mytheorem}\label{main2}
Let $G$ be a graph of order $n$ with no component of order less than $3$ and with all the bipartite components having both color classes of even order. Then:
$$
k(G)=\left\{\begin{array}{lll}
n,&\text{if}& n \not \equiv 2 \imod 4, \\
n+1,&\text{if}& n\equiv 2 \imod 4.
\end{array}
\right.
$$
\end{mytheorem}

\begin{mytheorem}\label{main3}
Let $G$ be a graph of order $n$ having neither component of order less than $3$ nor a $K_{1,2u+1}$ component for any integer $u\geq 1$. Then:
$$
\begin{array}{lll}
k(G)=n,& \text{if} &n \equiv 1\imod 2,\\
k(G)=n+1,& \text{if} &n\equiv 2 \imod 4,\\
k(G)\leq n+1, & \text{if} &n\equiv 0 \imod 4.
\end{array}
$$
\end{mytheorem}

\section{Proofs of the main results}

We will start this section with two simple observations.

\begin{myobservation}\label{pierwiastek}
If $\gr$ is an Abelian group of odd order $n$, then for any element $a\in\gr$ there exists unique element $b\in \gr$ such that $2b=a$.
\end{myobservation}

\noindent\textbf{Proof.}
Let us define the function $\lambda:\gr\rightarrow\gr$, $\lambda(g)=2g$. If for some $g_1,g_2\in\gr$ we have $\lambda(g_1)=\lambda(g_2)$, then it follows that $2(g_1-g_2)=0$ and consequently $g_1-g_2=0$, as there is no involution in $\gr$. Thus $\lambda$ is bijection.\qed\\

For sake of simplicity, for any element $a\in\gr$, we are going to use the notation $a/2$ for the element $b\in\gr$ satisfying $2b=a$. Let $\mu(\gr)$ be the number of elements $g\in \gr$ such that there exists $\frac{g}{2}\in\gr$.

\begin{myobservation}\label{parzysta}
Let $\gr$ be an Abelian group of order $n=2^{\alpha}(2\beta+1)$, for some integers $\alpha>0$ and $\beta$ and let $A$ be an Abelian group of order $2\beta+1$. If $\gr \cong  \zet_{2^{\alpha}}\times A$, then $\mu(\gr)=\frac{|\gr|}{2}$. Moreover if $i\neq 0$ is such an element of $\gr$ that $2i=0$, then there exists $\frac{i}{2} \in \gr$ if and only if $\alpha>1$.
\end{myobservation}

\noindent\textbf{Proof.}
Notice that for any  $a \in \zet_{2^\gamma}$, $a=2h+1$ for some $h\in \zet_{2^\gamma}$, there is no such element $a'\in \zet_{2^\gamma}$ that $2a'=a$.
Thus if $g=(2h+1,a)$, then there does not exist $\frac{g}{2} \in \gr$. On the other hand, if $g=(2h,a)$ for $a \in A$, then by Observation~\ref{pierwiastek} there exists $\frac{a}{2} \in A$ and thus $\frac{g}{2}=(h,\frac{a}{2})$. As it can be easily seen, there are exactly $2^{\alpha -1}(2k+1)$ elements $g\in \gr$ such that there exists $\frac{g}{2}$. Notice that if $i\neq 0$ is such element of $\gr$ that $2i=0$, then $i=(2^{\alpha-1},0)$. Thus there exists $\frac{i}{2} \in \gr$ if and only if $\alpha>1$.~\qed\\

Given any two vertices $x_1$ and $x_2$ belonging to the same connected component of $G$, there exist walks from $x_1$ to $x_2$. Some of them may consist of even number of vertices (some of them being repetitions). We are going to call them \textit{even walks}. The walks with odd number of vertices will be called \textit{odd walks}. We will always choose the shortest even or the shortest odd walk from $x_1$ to $x_2$.

We start with $0$ on all the edges of $G$. Then, in every step we will choose $x_1$ and $x_2$ and add some labels to all the edges of chosen walk from $x_1$ to $x_2$. To be more specific, we will add some element $a$ of the group to the labels of all the edges having odd position on the walk (starting from $x_1$) and $-a$ to the labels of all the edges having even position. It is possible that some labels will be modified more than once, as the walk does not need to be a path. We will denote such situation with $\phi_e(x_1,x_2)=a$ if we label the shortest even walk and $\phi_o(x_1,x_2)=a$ if we label the shortest odd walk. Observe that putting $\phi_e(x_1,x_2)=a$ means adding $a$ to the weighted degrees of both vertices, while $\phi_o(x_1,x_2)=a$ means adding $a$ to the weighted degree of $x_1$ and $-a$ to the weighted degree of $x_2$. In both cases the operation does not change the weighted degree of any other vertex of the walk. Note that if some component $G_1$ of $G$ is not bipartite, then for any vertices $x_1,x_2\in G_1$ there exist both even and odd walks.

Next lemma gives us the sufficient conditions for a group $\gr$ to allow some $\gr$-irregular labeling of given graph $G$.

\begin{mytheorem}\label{inwolucja}
Let $G$ be a disconnected graph of order $n$ with all the bipartite components having both color classes of even order and with no component of order less than $3$. Let $s=n+1$ if $n\equiv 2\imod 4$ and $s=n$ otherwise. Then for every integer $t\geq s$ there exists a $\gr$-irregular labeling for every Abelian group $\gr$ of order $t$ with at most one element $i\in\gr$ such that $i\neq 0=2i$. 
\end{mytheorem}

\noindent\textbf{Proof.} Let us denote the connected components of $G$ with $G_1,G_2,\dots,G_r$. Let $G_1,G_2,\ldots,G_p$ be the components of odd order, $G_{p+1},G_{p+2},\ldots,G_{q}$ the non-bipartite components of even order and $G_{q+1},G_{q+2},\ldots,G_{r}$ - the bipartite components.

As there is at most one involution in $\gr$, we can be sure that $\gr\equiv\zet_{2^\alpha}\times A$ where $A$ is an Abelian group of order $2\beta+1$ for some integers $\alpha,\beta\geq 0$. Let the elements of $\gr$ be $g_0=0, g_1, g_2,\ldots, g_{t-1}$. We can thus represent the elements of $\gr$ by the ordered pairs $(z,a)$, where $z\in \zet_{2^\alpha}$ and $a\in A$. We will use $a$ instead of $(0,a)$ when $\alpha=0$ (i.e. when $|\gr|=2\beta+1$).

Notice that $t\geq n+1$ for $n \equiv 2 \imod 4$ and $t\geq n$ otherwise. Thus the number of elements $g\in\gr$ for which $g/2$ exists is $\mu(\gr)\geq \frac{|\gr|}{2}\geq \frac{n}{2}>\frac{n}{3}\geq r\geq p$. Notice that if there exists $\frac{g}{2}\in \gr$, then exists $\frac{-g}{2}\in \gr$ and $\frac{g}{2}\neq \frac{-g}{2}$ if $g\neq 0$. Therefore we can pick $k=\lfloor\frac{p}{2}\rfloor$ elements, say $g_1, g_2,\ldots, g_k$, such that $g_j\neq 0$, $g_j\neq -g_j$ and there exists $\frac{g_j}{2}\in \gr$ for $j=1,2,\ldots,k$. We choose any vertex $x_j\in G_j$ for $j=1,\dots,2k-2$ and we put

$$
\phi_e(x_j,x_j)=\left\{
\begin{array}{lll}
\frac{g_{(j+1)/2}}{2},&\text{if}&j\equiv 1\imod 2,\\
\frac{-g_{j/2}}{2},&\text{if}&j\equiv 0\imod 2.
\end{array}
\right.
$$

Now we consider two cases, depending on the parity of $n$. If $n$ is odd, then also $p$ is odd. In such a situation we choose any vertex $x_j\in G_j$ for $j=2k-1,2k=p-1$ and we put

$$
\begin{array}{l}
\phi_e(x_{2k-1},x_{2k-1})=\frac{g_k}{2}\\
\phi_e(x_{2k},x_{2k})=\frac{-g_k}{2}
\end{array}
$$
and we leave the vertex $x_p$ with weighted degree $0$. After that we are left with even number $n-p$ of vertices that can be joined with even walks into couples $(x^1_j,x^2_j)$, $j=1,\dots,(n-p)/2$. On the other hand there are at least $(n-p)/2$ disjoint pairs $(g_j,-g_j)$, $g_j\neq 0$ of unused elements of $\gr$: $(|\gr|-p)/2$ when $|\gr|$ is odd and $(|\gr|-p-1)/2$ otherwise. We put
$$
\phi_o(x^1_j,x^2_j)=g_j
$$
for every $j=1,\dots,(n-p)/2$ and we are done.

If $n$ is even, then also $p$ is even. Our next steps depend on the value of $t$. If $t\geq n+1$, then we proceed as in the previous case, except the fact that $p=2k$ and none of the vertex weights in odd components is equal to $0$. Observe that also in this case after weighting one vertex in every odd component we still have at least $(n-p)/2$ disjoint pairs $(g_j,-g_j)$, $g_j\neq 0$ of unused elements of $\gr$ and $n-p$ vertices to label. If $t=n$, it means that $n\equiv 0\imod 4$ and so there exists $i/2\in\gr$, where $i$ is the only involution of $\gr$. So we choose any vertex $x_{2k-1}\in G_{2k-1}$, we put
$$
\phi_e(x_{2k-1},x_{2k-1})=\frac{i}{2}
$$
and we leave the vertex $x_p$ with weighted degree $0$. After that we still have exactly $(n-p)/2$ disjoint pairs $(g_j,-g_j)$, $g_j\neq 0$ of unused elements of $\gr$ and $n-p$ vertices to label, so we proceed as in the previous cases and we can complete the labeling.
\qed\\

In particular, the following corollary is true, as in every cyclic group $\zet_{k}$ there is exactly one involution $i=k/2$ if $k$ is even and no involution if $k$ is odd.

\begin{mycorollary}\label{cykliczna}
Let $G$ be a disconnected graph of order $n$ with all the bipartite components having both color classes of even order and with no component of order less than $3$. Let $s=n+1$ if $n\equiv 2\imod 4$ and $s=n$ otherwise. Then for every integer $t\geq s$ there exists a $\zet_t$-irregular labeling. 
\end{mycorollary}

Next lemma allows to find group irregular labelings for another class of graphs.

\begin{mylemma}\label{lemZSP}
Let $G$ be a graph of order $n$ having neither component of order less than $3$ nor a $K_{1,2u+1}$ component for any integer $u\geq 1$. Then for every odd integer $t\geq n$ there exists a $\zet_t$-irregular labeling.
\end{mylemma}  

\noindent\textbf{Proof.}
We are going to use the following theorem, proved in \cite{ref_KapLevRod}.

\begin{mytheorem}[\cite{ref_KapLevRod}]\label{thmZSP}
Let $n=r_1+r_2+\dots+r_q$ be a partition of the positive integer $t$, where $r_i\geq 2$ for $i=1,2,\dots,q$. Let $A=\{1,2,\dots,t\}$.
Then the set $A$ can be partitioned into pairwise disjoint subsets $A_1, A_2, \dots, A_q$ such that for every $1\leq i\leq q$, $|A_i|=r_i$ with
$\sum_{a\in A_i}{a}\equiv 0\imod {t+1}$ if $t$ is even and $\sum_{a\in A_i}{a}\equiv 0\imod t$ if $t$ is odd.
\end{mytheorem}

We are going to divide the vertices of $G$ into triples and pairs. Let $p_1$ be the number of bipartite components of $G$ with both color classes odd, $p_2$ with both classes even and $p_3$ with one class odd and one even. Let $p_4$ be the number of remaining components of odd order and $p_5$ - the number of remaining components of even order. The number of triples equals to $2p_1+p_3+p_4$. The remaining vertices form the pairs.

Now we are going to partition the elements of $\zet_{t}$. We are interested only in the case when $t$ is odd. The authors show the construction of $2l+1$ triples $B_1,B_2,\dots,B_{2l+1}$ and $m$ pairs $C_1,C_2,\dots,C_m$ where $l=\lfloor(2p_1+p_3+p_4)/2\rfloor$ and $m=(t-6l-3)/2$. We will denote the elements of triples and pairs with $B_j=(a_j,b_j,c_j)$ and $C_j=(d_j,-d_j)$. Each of those triples and pairs sums up to $0\imod t$. It is easy to observe that for a given element $g\in\zet_t$, $g\neq -g$, either $(g,-g)=C_j$ for some $j$ or $(g,-g,0)=B_j$ for some $j$ (there exists exactly one such triple) or $g$ and $-g$ belong to two distinct triples. If $n$ is even, then we remove $t$ (i.e. $0$) from the triple $B_j$ to which it belongs in order to obtain additional pair $C_{m+1}$ (observe that after removing $0$, the number of remaining elements of $\zet_t$ is still enough to distinguish all the vertices of $G$ as in this case $t\geq n+1$).

Let us start the labeling. We are numbering the pairs and triples consecutively, starting with $1$ for both lists.

Given any bipartite component $G$ with both color classes even, we divide the vertices of every color class into pairs $(x_j^1,x_j^2)$, putting
$$
\phi_o(x_j^1,x_j^2)=d_j
$$
for every such pair. We proceed in similar way in the case of all the non-bipartite components of even order, coupling the vertices of every such component in any way.

If both color classes of a bipartite component are of odd order, then they both have at least $3$ vertices. We choose three of them, denoted with $x_j$, $y_j$ and $z_{j+1}$, in one class and another three, $x_{j+1}$, $y_{j+1}$ and $z_j$, in another one and we put
$$
\begin{array}{l}
\phi_e(x_j,z_j)=b_j,\\
\phi_e(y_j,z_j)=c_j,\\
\phi_e(x_{j+1},z_{j+1})=b_{j+1},\\
\phi_e(y_{j+1},z_{j+1})=c_{j+1}.
\end{array}
$$
We proceed with the remaining vertices of this kind of components as in the case when both color classes are even.

In the case of all the components of odd order we choose three vertices $x_j$, $y_j$ and $z_j$ (in the case of bipartite component $z_j$ belongs to the odd color class and to other vertices to the even one). We put
$$
\begin{array}{l}
\phi_e(x_j,z_j)=b_j,\\
\phi_e(y_j,z_j)=c_j.
\end{array}
$$
Observe that the numbers of triples and pairs of elements of $\zet_t$ are at least equal to the numbers of triples and pairs of vertices of $G$. Thus the labeling defined above is $zet_t$-irregular. Indeed, in the $j^{th}$ triple of vertices the weights are equal to $w(x_j)=b_j$, $w(y_j)=c_j$ and $w(z_j)=-a_j$ and in the $j^{th}$ pair we have $w(x_j^1)=d_j$ and $w(x_j^2)=-d_j$.
\qed\\

The main results easily follow from the above theorems and observations. Let $\gr$ be a group of order $t$, where $t=n$ if $n$ is odd and $t=n+1$ otherwise. In both cases $t$ is odd and there is no element $i\in \gr$ such that $i \neq 0=2i$. Thus the Theorem \ref{main} easily follows from the Theorem~\ref{inwolucja} and Lemma~\ref{lemma_below}. The Theorem~\ref{main2} follows from Corollary \ref{cykliczna} and Lemma~\ref{lemma_below}.  Finally, the Theorem~\ref{main3} is consequence of Lemma \ref{lemZSP} and Lemma~\ref{lemma_below}.

\end {document}